\documentclass[10pt,amssymb,amstex]{amsart}
\numberwithin{equation}{section}

\newtheorem{thm}{Theorem}[section]

\newtheorem{defin}[thm]{Definition}

\begin{document}

\author{Ravshan Ashurov}
\address{Ashurov R: Institute of Mathematics of Uzbekistan,
Tashkent, Uzbekistan.}
\address{AU Engineering School Akfa University 264, Milliy Bog Str. Tashkent 111221 Uzbekistan.}
\email{ashurovr@gmail.com}

\author{Yusuf Fayziev}
\address{Fayziev Yu: National University of Uzbekistan,  Tashkent, Uzbekistan.}
\email{fayziev.yusuf@mail.ru}

\author{Muattar Khudoykulova}
\address{Khudoykulova M: National University of Uzbekistan,  Tashkent, Uzbekistan.}
\email{muattarxudoyqulova2000@gmail.com}

\small

\title[Forward and inverse problems for subdiffusion equation...] {Forward and inverse problems for subdiffusion equation with  time-dependent coefficients}

\begin{abstract}
In this paper, we consider forward and inverse problems for subdiffusion equations with time-dependent coefficients. The fractional derivative is taken in the sense of Riemann-Liouville. Using the classical Fourier method, the theorem of the uniqueness and existence of forward and inverse problems for determining the right-hand side of the equation are proved.

\end{abstract}

\maketitle

\section{Main results}

In modern science and engineering researchers frequently apply fractional order differential equations to model  dynamics of various processes arising in different fields; see, for example, \cite{Hilfer,Mainardi,MetzlerKlafter00,West} in physics, \cite{MachadoLopes,SGM} in finance,   \cite{BMR} in hydrology, \cite{Magin} in cell biology,  among others. In the last few decades  several books,  devoted to fractional order differential equations and their various applications, have been published (see e.g.
\cite{SKM,Podlubny,KST,Handbook,Umarov,UHK_book_2018}).
However, the state of fractional differential equations with time-dependent coefficients is still far from completeness, even for subdiffusion equations.
This is due to many reasons. For example, if the coefficients of the elliptic operator in the subdiffusion equation depend on $t$, then, as a rule, it is not possible to apply the separation of variables method to the corresponding initial-boundary value problem. Another reason is that the fractional equation with respect to the variable $t$, obtained after separation of variables, cannot be solved with using known methods.

Initial-boundary value and the Cauchy problems for a subdiffution equations have been considered by a
number of authors using different methods. Many works are devoted to the case of one space variable $x\in\mathbb{R}$ and the subdiffusion equation with the "elliptic part" $u_{xx}$ (see, for
example, monograph of A. V. Pskhu \cite{Pskh} and literature therein).

 In multidimensional case ($x\in \mathbb{R}^N$)
instead of the differential expression $u_{xx}$ it has been
considered elliptic operators with constant coefficients in books \cite{KST}, \cite{Gor} (and literature therein) and the survey paper \cite{AK} (and literature therein).

Authors of \cite{KST}, \cite{Gor} and \cite{AK}  investigated the Cauchy type
problems applying either the Laplace transform or the Fourier
transform.

In his recent paper \cite{Pskh} A. V. Pskhu considered
an initial-boundary value problem for subdiffusion equation with
the Laplace operator and domain $\Omega$ - a multidimensional
rectangular region. The author succeeded to construct the Green's
function.

In an arbitrary domain $\Omega$ initial-boundary value problems
for subdiffusion equations (the fractional part of the equation is
a multi-term and initial conditions are non-local) with the Caputo
derivatives has been investigated by M. Ruzhansky et al \cite{RuTT}. The authors proved the existence and uniqueness of the generalized solution to the problem.The existance and uniqueness of the classical solution for this general case was proved in \cite{AOLob}.

There many other papers, where it is proved the existence and  uniqueness of the generalized and the weak solutions, see the survey paper\cite{Rico}  and literature therein.

	As for subdiffusion equations with time-dependent coefficients, we can cite the work of Kubica and Yamamoto \cite{Kubica} (see also \cite{RicoBOSHQA}), where the existence of a weak solution for an initial-boundary value problem with an elliptic operator, whose coefficients depend on $t$ is investigated. Since the Fourier method is not applicable, the existence of a solution was proved by the method of a priori estimates.

In the present paper, we study the existence of a strong solution to the initial-boundary value problem for a certain class of subdiffusion equations with coefficients depending on $t$. Also, for such equations, we study the inverse problem of determining the right-hand side of the equation.

Let us pass to the exact formulations of the results obtained.

Let $H$ be a separable Hilbert space with the scalar product
$(\cdot, \cdot)$ and the norm $||\cdot||$ and  $A: H\rightarrow H$
be an arbitrary unbounded positive selfadjoint operator in $H$, with the domain of definition $D(A)$.
Suppose that $A$ has a complete in $H$ system of orthonormal
eigenfunctions $\{v_k\}$ and a countable set of nonnegative
eigenvalues $\lambda_k$. It is convenient to assume that the
eigenvalues do not decrease as their number increases, i.e.
$0<\lambda_1\leq\lambda_2 \cdot\cdot\cdot\rightarrow +\infty$.

Let $\alpha\in(0,1) $ and $C((a,b); H)$ stand
for a set of continuously differentiable  functions $h(t)$  of $t\in (a,b)$ with values in $H$. Using the definitions of a strong integral and a strong
derivative, fractional analogues of derivatives can
be determined for vector-valued functions (or simply functions)
$h: \mathbb{R}_+\rightarrow H$, while the well-known formulae and
properties are preserved (see, for example, \cite{Liz}). Recall
that the Riemann-Liouvelle fractional derivative (and integral) of order $ \alpha $ ($\sigma<0$) of the
function $ h (t) $ defined on $ [0, \infty) $ has the form
\begin{equation}\label{def0}
\partial_t^\alpha h(t)=\frac{d}{dt} J^{\alpha-1} h(t),\,\, J^\sigma h(t) = \frac{1}{\Gamma
	(-\sigma)}\int\limits_0^t\frac{h(\xi)}{(t-\xi)^{\sigma+1}} d\xi,
\quad t>0,
\end{equation}
provided the right-hand side exists. Here $\Gamma(\sigma)$ is
Euler's gamma function.

Consider the following problem
\begin{equation}\label{prob1}
\left\{
\begin{aligned}
&\partial_t^\alpha u(t) + t^{\beta}Au(t) = t^{\mu}f,\quad 0<t\leq T; \\
&\lim\limits_{t\rightarrow +0}J_t^{\alpha-1} u(t) =
\varphi,
\end{aligned}
\right.
\end{equation}
where $f,\varphi\in H$ and $\beta>-\alpha$, $\mu>-1$. This problem are also called \emph{the forward problem}.

\begin{defin}\label{def}

 A function \,  $u(t)\in C((0,\infty); H)$  \, with the properties
$J^{\alpha-1}_t u(t)\in C([0,\infty); H)$, \, $\partial_t^\alpha u(t)$, \, \, $Au(t)\in C((0,\infty); H)$ and
satisfying conditions (\ref{prob1})  is called \textbf{the
solution} of problem (\ref{prob1}).
\end{defin}

To formulate our result on the forward problem we need to remind the definition of the Mittag-Leffler function with three
parameters $
E_{\alpha, m, l}(z)$, where $m$ and $l$ are an arbitrary complex numbers (see e.g. \cite{KST}). We have
\begin{equation}\label{ml}
	E_{\alpha, m, l}(z)= \sum\limits_{k=0}^\infty {c_k}{z^k},
\end{equation}
where
\begin{equation}\label{m1}
	{c_0}=1, \quad
	{c_k}=\prod_{j=0}^{k-1} \frac{\Gamma(\alpha(jm+l)+1)}{\Gamma(\alpha(jm+l+1)+1)}.
\end{equation}

\begin{thm}\label{th1}
 Let  $\varphi$ and $f \in H$. Then problem (\ref{prob1}) has a unique solution and this solution has the form
$$
u(t)=\sum\limits_{k=1}^\infty \left[\frac{\varphi_{k}}{\Gamma(\alpha)} t^{\alpha-1} E_{\alpha,1+\frac{\beta}{\alpha},1+\frac{\beta-1}{\alpha}} (-\lambda_k
t^{\alpha+\beta})\right.
$$
\begin{equation}\label{ux1}
\left.+\frac{\Gamma(\mu+1){f_{k}}t^{\alpha+\mu}}{\Gamma(\mu+\alpha+1)} E_{\alpha,1+\frac{\beta}{\alpha},1+\frac{\beta+\mu}{\alpha}}(-\lambda_{k}t^{\alpha+\beta})\right]v_k,
\end{equation}
where $f_k$,  $\varphi_{k}$ are the Fourier coefficients of the functions $f$ and $\varphi$ respectively.

\end{thm}

Next, we will study the inverse problem of determining the right-hand side of equation (\ref{prob1}). In this case, as a over-determination condition, we set the value of the solution in final time:
\begin{equation}\label{prob4}
    u(T)=\Psi.
\end{equation}
The inverse problems of determining the right-hand side (the heat
source density) of various subdiffusion equations have been
considered by a number of authors (see, e.g. \cite{12,25,KiraneM,KiraneST,13}
and the bibliography therein). An overview of works in this direction published up to 2019 can be found in \cite{YLY}, p. 411-431.

Let us dwell on some works published after 2019.

 The paper \cite{AODif} is devoted to the inverse
problem for the subdiffusion equation with Riemann-Liouville
derivatives, the elliptical part of which has the most general form.

In \cite{AF1}  the authors
considered an inverse problem for the simultaneous determination
of the order of the Riemann-Liouville  fractional derivative and
the source function in the subdiffusion equations. Using the
classical Fourier method, the authors proved the uniqueness and
existence theorem for this inverse problem.

In works \cite{AF4} and \cite{AFKaraganda}, the authors studied, along with other problems, the inverse problem of determining the right-hand side of the subdiffusion equations for non-local problems in time.

Authors of the work \cite{ASH} considered the inverse problem of determining the time-dependent source function for Schrödinger type equations of fractional order. As a over-determination condition for solving this inverse problem, the authors have taken condition $B[u(\cdot, t)]=\psi(t)$, where $B$ is an arbitrary linear bounded functional. The existence and uniqueness theorem was proved.

In \cite{AF2,AF3}, the authors investigated the inverse
problem of determining the order of the fractional derivative in
the subdiffusion equation and in the wave equation, respectively.

\begin{defin}\label{defIN}
     A pair $\{u(t), f \}$ of function $u(t)\in C((0,\infty); H)$ and $f \in H$ with the properties $J^{\alpha-1}_t u(t)\in C([0,\infty); H)$, \,
$\partial_t^\alpha u(t), Au(t)\in C((0,\infty); H)$ and satisfying
conditions (\ref{prob1}), (\ref{prob4})  is called \textbf{the
 solution} of the inverse  problem (\ref{prob1}), (\ref{prob4}).

\end{defin}

We have
\begin{thm}\label{th4}
     Let  $\varphi, \Psi \in D(A)$. Then the inverse
problem (\ref{prob1}), (\ref{prob4}) has a unique solution
$\{u(t), f \}$ and this solution has the form
$$
u(t)=\sum\limits_{k=1}^\infty \left[\frac{\varphi_{k}}{\Gamma(\alpha)}t^{\alpha-1} E_{\alpha,1+\frac{\beta}{\alpha},1+\frac{\beta-1}{\alpha}} (-\lambda_k
t^{\alpha+\beta})\right.
$$
\begin{equation}\label{ux}
\left.+\frac{\Gamma(\mu+1){f_{k}}}{\Gamma(\mu+\alpha+1)} t^{\alpha+\mu} E_{\alpha,1+\frac{\beta}{\alpha},1+\frac{\beta+\mu}{\alpha}}(-\lambda_{k}t^{\alpha+\beta})\right]
v_k,
\end{equation}
where the Fourier coefficients $f_k$ have the form
$$
f_k= \frac{\Psi_k \Gamma(\mu+\alpha+1)}{\Gamma(\mu+1)T^{\alpha+\mu} E_{\alpha,1+\frac{\beta}{\alpha},1+\frac{\beta+\mu}{\alpha}}(-\lambda_{k}T^{\alpha+\beta})}
$$
\begin{equation}\label{fk}
-\frac{\varphi_{k}T^{\alpha-1}E_{\alpha,1+\frac{\beta}{\alpha},\frac{\beta}{\alpha}} (-\lambda_k
T^{\alpha+\beta})\Gamma(\mu+\alpha+1)}{\Gamma(\alpha)
\Gamma(\mu+1)T^{\alpha+\mu} E_{\alpha,1+\frac{\beta}{\alpha},1+\frac{\beta+\mu}{\alpha}}(-\lambda_{k}T^{\alpha+\beta})}
\end{equation}
and
\begin{equation}\label{uff}
f=\sum\limits_{k=1}^\infty f_k v_k.
\end{equation}

\end{thm}

\section{Preliminaries}

Let us introduce the class of "smooth" functions with the help of the degree of operator $A$.
We set for an arbitrary real number $ \tau $
$$
A^\tau h= \sum\limits_{k=1}^\infty \lambda_k^\tau h_k v_k,
$$
where $h_k$ is the Fourier coefficients of a function $h\in H$:
$h_k=(h, v_k)$. Obviously, the domain of this operator has the
form
$$
D(A^\tau)=\{h\in H:  \sum\limits_{k=1}^\infty \lambda_k^{2\tau}
|h_k|^2 < \infty\}.
$$
For elements of $D(A^\tau)$ we introduce the norm
\[
||h||^2_\tau=\sum\limits_{k=1}^\infty \lambda_k^{2\tau} |h_k|^2 =
||A^\tau h||^2,
\]
and together with this norm $D(A^\tau)$ turns into a Hilbert
space.

Let us remind some properties of the Mittag-Leffler function (see,\cite{30},  p.9 and 16 )

\label{pro2} \textbf{Proportion 1.} Let $\alpha \in (0,1)$,\, $m>0, \,\, t\geq 0$. Then
\begin{equation}\label{ML1}
E_{\alpha, m, m-\frac{1}{\alpha}}(-t)\leq \frac{1}{\left(1+\frac{\Gamma(1+\alpha m)}{\Gamma(1+\alpha(1+m))}t\right)^{1+\frac{1}{m}}}.
\end{equation}

\label{pro2} \textbf{Proportion 2.} Let $\alpha \in (0,1)$,\, $m>0,\, l>m-\frac{1}{\alpha}, \,\, t\geq 0$. Then
\begin{equation}\label{ML2}
\frac{1}{1+\frac{\Gamma(1+\alpha(l-m))}{\Gamma(1+\alpha(l-m+1))}t}\leq E_{\alpha, m, l}(-t)\leq \frac{1}{1+\frac{\Gamma(1+\alpha l)}{\Gamma(1+\alpha(1+l))}t}.
\end{equation}

\section{Proof of Theorem \ref{th1}}

Let us prove the existence and uniqueness
theorem for solutions of problem (\ref{prob1}).

{\em Proof.} \textbf{Existence}. We use the Fourier method to solve problem (\ref{prob1}) and look for a solution in the form

 \[
	u(t) = \sum\limits_{k=1}^\infty T_k(t) v_k,
	\]
where $T_k(t)$, $k\geq 1$, are  solutions of the Cauchy problem:
	
 \begin{equation}\label{probT}
\left\{
 \begin{aligned}
&\partial_t^\alpha T_k(t) + \lambda_k t^{\beta}T_k(t) =t^{\mu}f_k,\quad t>0; \\
&\lim\limits_{t\rightarrow +0}J_t^{\alpha-1} T_k(t) =
\varphi_k.
\end{aligned}
\right.
\end{equation}
Problem (\ref{probT}) has a unique solution  (see, e.g. \cite{KST}, p. 246), and the solution has the form:
 \begin{equation}\label{ST}
T_k(t)= \frac{\varphi_{k}}{\Gamma(\alpha)} t^{\alpha-1} E_{\alpha,1+\frac{\beta}{\alpha},1+\frac{\beta-1}{\alpha}} (-\lambda_k
t^{\alpha+\beta})+\frac{\Gamma(\mu+1){f_{k}}}{\Gamma(\mu+\alpha+1)} t^{\alpha+\mu} E_{\alpha,1+\frac{\beta}{\alpha},1+\frac{\beta+\mu}{\alpha}}(-\lambda_{k}t^{\alpha+\beta}).
\end{equation}

Therefore, series (\ref{ux1}) is a formal solution to problem (\ref{prob1}). Let us now show that this series and the series obtained after differentiation converge.

Let $S_n(t)$ be the partial sum of series (\ref{ux1}). Due to the Parseval
equality we may write
$$
||S_n(t)||^2 = \sum\limits_{k=1}^n
\left|\frac{\varphi_{k}}{\Gamma(\alpha)}t^{\alpha-1} E_{\alpha,1+\frac{\beta}{\alpha},1+\frac{\beta-1}{\alpha}} (-\lambda_k
t^{\alpha+\beta})\right.
$$
$$
\left.+\frac{\Gamma(\mu+1){f_{k}}}{\Gamma(\mu+\alpha+1)} t^{\alpha+\mu} E_{\alpha,1+\frac{\beta}{\alpha},1+\frac{\beta+\mu}{\alpha}}(-\lambda_{k}t^{\alpha+\beta})\right|^2.
$$
Then, we have
$$
||S_n(t)||^2\leq  2\sum\limits_{k=1}^n
|\frac{\varphi_{k}}{\Gamma(\alpha)}t^{\alpha-1} E_{\alpha,1+\frac{\beta}{\alpha},1+\frac{\beta-1}{\alpha}} (-\lambda_k
t^{\alpha+\beta})|^2
$$
$$
+2\sum\limits_{k=1}^n |\frac{\Gamma(\mu+1){f_{k}}}{\Gamma(\mu+\alpha+1)} t^{\alpha+\mu} E_{\alpha,1+\frac{\beta}{\alpha},1+\frac{\beta+\mu}{\alpha}}(-\lambda_{k}t^{\alpha+\beta})|^2=
2  S_n^1+2 S_n^2,
$$
where
$$
S_n^1=\sum\limits_{k=1}^n
|\frac{\varphi_{k}}{\Gamma(\alpha)}t^{\alpha-1} E_{\alpha,1+\frac{\beta}{\alpha},1+\frac{\beta-1}{\alpha}} (-\lambda_k
t^{\alpha+\beta})|^2,
$$
and
$$
S_n^2=\sum\limits_{k=1}^n \left|\frac{\Gamma(\mu+1){f_{k}}}{\Gamma(\mu+\alpha+1)} t^{\alpha+\mu} E_{\alpha,1+\frac{\beta}{\alpha},1+\frac{\beta+\mu}{\alpha}}(-\lambda_{k}t^{\alpha+\beta})\right|^2.
$$
Apply inequality (\ref{ML1}) to get
$$
S_n^1\leq C\sum\limits_{k=1}^n
|\varphi_{k}|^2 t^{2(\alpha-1)} \left|\frac{1}{1+\frac{\Gamma(\alpha+\beta+1)}{\Gamma(2\alpha+\beta+1)}\lambda_k
t^{\alpha+\beta}}\right|^{\frac{4\alpha}{\alpha+\beta}}\leq \frac{C_1} {\lambda_1^{\frac{4\alpha}{\alpha+\beta}}t^{2(\alpha+1)}}\sum\limits_{k=1}^n
|\varphi_{k}|^2,
$$
and
$$
S_{n}^2\leq C\left(\frac{\Gamma(\mu+1)}{\Gamma(\mu+\alpha+1)}\right)^2t^{2(\alpha+\mu)}\sum\limits_{k=1}^n |f_k|^2  \left|\frac{1}{1+\frac{\Gamma(1+\alpha+\beta+\mu)}{\Gamma(1+2\alpha+\beta+\mu)}\lambda_k t^{\alpha+\beta}}\right|^2
$$
$$
\leq C\left(\frac{\Gamma(\mu+1)\Gamma(1+2\alpha+\beta+\mu)}{\Gamma(\mu+\alpha+1)\Gamma(1+\alpha+\beta+\mu)}\right)^2\frac{1}{\lambda_1^2 t^{2(\beta-\mu)}}\sum\limits_{k=1}^n
|f_k|^2=
\frac{C_1}{t^{2(\beta-\mu)}} \sum\limits_{k=1}^n|f_k|^2.
$$
Since $\varphi, f \in H$  then sum (\ref{ux1}) converges and
$u(t)\in C((0,\infty); H)$.

Now let us estimate $Au(t)$. One has
$$
AS_n(t) = \sum\limits_{k=1}^n
\left[\frac{\varphi_{k}}{\Gamma(\alpha)}t^{\alpha-1} E_{\alpha,1+\frac{\beta}{\alpha},1+\frac{\beta-1}{\alpha}} (-\lambda_k
t^{\alpha+\beta})\right.
$$
\begin{equation}
\left.+\frac{\Gamma(\mu+1){f_{k}}}{\Gamma(\mu+\alpha+1)} t^{\alpha+\mu} E_{\alpha,1+\frac{\beta}{\alpha},1+\frac{\beta+\mu}{\alpha}}(-\lambda_{k}t^{\alpha+\beta})\right]\lambda_kv_k.
\end{equation}
Apply the Parseval equality to obtain
$$
||AS_n(t)||^2 = \sum\limits_{k=1}^n\lambda_k^2
\left|\frac{\varphi_{k}}{\Gamma(\alpha)}t^{\alpha-1} E_{\alpha,1+\frac{\beta}{\alpha},1+\frac{\beta-1}{\alpha}} (-\lambda_k
t^{\alpha+\beta})\right.
$$
$$\left.+\frac{\Gamma(\mu+1){f_{k}}}{\Gamma(\mu+\alpha+1)} t^{\alpha+\mu} E_{\alpha,1+\frac{\beta}{\alpha},1+\frac{\beta+\mu}{\alpha}}(-\lambda_{k}t^{\alpha+\beta})\right|^2.
$$
Then, we have
$$
||AS_n(t)||^2\leq
C\sum\limits_{k=1}^n\lambda_k^2
\left|\frac{\varphi_{k}}{\Gamma(\alpha)}t^{\alpha-1} E_{\alpha,1+\frac{\beta}{\alpha},1+\frac{\beta-1}{\alpha}} (-\lambda_k
t^{\alpha+\beta})\right|^2
$$
$$
+C\sum\limits_{k=1}^n\lambda_k^2 \left|\frac{\Gamma(\mu+1){f_{k}}}{\Gamma(\mu+\alpha+1)} t^{\alpha+\mu} E_{\alpha,1+\frac{\beta}{\alpha},1+\frac{\beta+\mu}{\alpha}}(-\lambda_{k}t^{\alpha+\beta})\right|^2=
 AS_n^1+AS_n^2.
$$
By virtue of inequality (\ref{ML2}) we obtain
$$
AS_n^1=C\sum\limits_{k=1}^n\lambda_k^2
\left|\frac{\varphi_{k}}{\Gamma(\alpha)}t^{\alpha-1} E_{\alpha,1+\frac{\beta}{\alpha},1+\frac{\beta-1}{\alpha}} (-\lambda_k
t^{\alpha+\beta})\right|^2
$$
$$ \leq C_1 t^{2(\alpha-1)} \sum\limits_{k=1}^n\lambda_k^2
|\varphi_{k}|^2 \left|\frac{1}{\left(1+\frac{\Gamma(\alpha+\beta+1)}{\Gamma(2\alpha+\beta+1)}\lambda_k
t^{\alpha+\beta}\right)^{\frac{2\alpha+\beta}{\alpha+\beta}}}\right|^2\leq
$$
$$
\leq \frac{C_2}{t^{2(\alpha+\beta+1)}}\sum\limits_{k=1}^n\lambda_k^2\frac{|\varphi_{k}|^2}{\lambda_k^{2+\frac{2\alpha}{\alpha+\beta}}}=\frac{C_2}{\lambda_1^{\frac{2\alpha}{\alpha+\beta}}t^{2(\alpha+\beta+1)}}\sum_{k=1}^n|\varphi_{k}|^2,
$$
and
$$
AS_n^2\leq C\left(\frac{\Gamma(\mu+1)}{\Gamma(\mu+\alpha+1)}\right)^2t^{2(\alpha+\mu)}\sum\limits_{k=1}^n\lambda_k^2
|f_k|^2  \left|\frac{1}{1+\frac{\Gamma(1+\alpha+\beta+\mu)}{\Gamma(1+2\alpha+\beta+\mu)}\lambda_k t^{\alpha+\beta}}\right|^2
$$
$$
\leq \frac{C}{t^{2(\beta-\mu)}}\left(\frac{\Gamma(\mu+1)\Gamma(1+2\alpha+\beta+\mu)}{\Gamma(\mu+\alpha+1)\Gamma(1+\alpha+\beta+\mu)}\right)^2 \sum\limits_{k=1}^n\lambda_k^2
\frac{|f_k|^2}{\lambda_k^2}\leq \frac{C}{t^{2(\beta-\mu)}}\sum\limits_{k=1}^n|f_k|^2.
$$
Again, since $\varphi_k, f \in H$  one has $Au(t)\in C((0,\infty);
H)$.

Further, using equations (\ref{prob1}), we can write $\partial_t^\alpha
S_n(t)= - AS_n(t)+\sum\limits_{k=1}^n f_k(t) v_k$, $t>0$. Since we already have an estimate for the right side, it follows that the series (\ref{ux1}) can be differentiated term-by-term with respect to the variable $t$.

\textbf{Uniqueness.} The uniqueness of the solution can be proved by the standard technique based on completeness of the set of
eigenfunctions $\{v_k\}$ in $H$. Indeed, let $u(t)$ be a solution to the homogeneous
problem:
\begin{equation}\label{eq200}
\left\{
\begin{aligned}
&\partial_t^\alpha u(t) + t^{\beta}Au(t) = 0,\quad t>0; \\
&\lim\limits_{t\rightarrow +0}J_t^{\alpha-1} u(t) =0.
\end{aligned}
\right.
\end{equation}
We set $u_k(t)=(u(t),
v_k)$. Then, by virtue of equation (\ref{eq200}) and the
selfadjointness of the operator $A$,
\begin{equation}\label{eq201}
\partial_t^\alpha u_k(t)=(\partial_t^\alpha u(t), v_k)= -t^{\beta}(Au(t), v_k)=-t^{\beta}(u(t), Av_k)
\end{equation}
$$=-t^{\beta}(u(t), \lambda_k v_k) =-\lambda_k t^{\beta}(u(t), v_k) =-\lambda_k t^{\beta} u_k(t), \quad
t>0.
$$
Thus, we have the following problem
$$
\partial_t^\alpha u_k(t) +\lambda_kt^{\beta} u_k(t)=0,\quad t>0; \quad J^{\alpha-1}_t u_k(+0) =0.
$$

 This problem has the unique solution: $u_k(t)\equiv 0$ for all  $k$ (see, e.g. \cite{KST}).

Consequently, due to the completeness of the system of
eigenfunctions $ \{v_k \} $, we have $u(t) \equiv 0$, as required.

\section{Proof of Theorem \ref{th4}}

{\em Proof.} \textbf{Existence. } As proved above, if $f$ is known, then the solution to the forward problem has the form (\ref{ux}). Now, using the over-determination condition (\ref{prob4}), we find the unknown right-hand side $f$.

By virtue of the
completeness of system $\{v_k\}$ we obtain:
$$
 \frac{\varphi_{k}}{\Gamma(\alpha)}T^{\alpha-1} E_{\alpha,1+\frac{\beta}{\alpha},1+\frac{\beta-1}{\alpha}} (-\lambda_k
T^{\alpha+\beta})
$$
$$+\frac{\Gamma(\mu+1){f_{k}}}{\Gamma(\mu+\alpha+1)} T^{\alpha+\mu} E_{\alpha,1+\frac{\beta}{\alpha},1+\frac{\beta+\mu}{\alpha}}(-\lambda_{k}T^{\alpha+\beta})=\Psi_k.
$$
After simple calculations, we get
$$
f_k= \frac{\Psi_k \Gamma(\mu+\alpha+1)}{\Gamma(\mu+1)T^{\alpha+\mu} E_{\alpha,1+\frac{\beta}{\alpha},1+\frac{\beta+\mu}{\alpha}}(-\lambda_{k}T^{\alpha+\beta})}
$$
\begin{equation}
-\frac{\varphi_{k}T^{\alpha-1}E_{\alpha,1+\frac{\beta}{\alpha},1+\frac{\beta-1}{\alpha}} (-\lambda_k
T^{\alpha+\beta})\Gamma(\mu+\alpha+1)}{\Gamma(\alpha)\Gamma(\mu+1)T^{\alpha+\mu} E_{\alpha,1+\frac{\beta}{\alpha},1+\frac{\beta+\mu}{\alpha}}(-\lambda_{k}T^{\alpha+\beta})}\equiv f_{k,1}+f_{k,2}.
\end{equation}
With these Fourier coefficients we have the above formal series
(\ref{uff}) for the unknown function $f$: $f=\sum\limits_{k=1}^\infty
\bigg(f_{k,1}+ f_{k,2}\bigg)v_k$.

Let us prove convergence of series (\ref{uff}). If $F_n$ is the
partial sums of series (\ref{uff}), then by virtue of the Parseval
equality we may write
\begin{equation}\label{Fj}
||F_n||^2 = \sum\limits_{k=1}^n \bigg[f_{k,1}+f_{k,2}\bigg]^2\leq 2
\sum\limits_{k=1}^n f_{k,1}^2+ 2\sum\limits_{k=1}^n f_{k,2}^2
\equiv 2 I_{1,n}+2 I_{2,n}.
\end{equation}
Then for $I_{1,n}$ one has
\[
I_{1,n} \leq  \sum\limits_{k=1}^n
{\left|\frac{\Psi_k \Gamma(\mu+\alpha+1)}{\Gamma(\mu+1)T^{\alpha+\mu} E_{\alpha,1+\frac{\beta}{\alpha},1+\frac{\beta+\mu}{\alpha}}(-\lambda_{k}T^{\alpha+\beta})}\right|^2}.
\]
Using (\ref{ML2}) we get
\[
I_{1,n} \leq C\sum\limits_{k=1}^n \frac{\lambda_k^2
|\Psi_k|^2}{\left( T^{\alpha+\mu}\frac{1}{T^{\alpha+\beta}}\right)^2} \leq C T^{2(\beta-\mu)}\sum\limits_{k=1}^n
\lambda_k^2 |\Psi_k|^2.
\]
Since $|E_{\alpha,m,l}(-\lambda_k T^\rho)|\leq 1$ we have
$$I_{2,n}\leq \sum\limits_{k=1}^n \left|\frac{\varphi_k T^{\alpha-1}E_{\alpha,1+\frac{\beta}{\alpha},1+\frac{\beta-1}{\alpha}} (-\lambda_k T^{\alpha+\beta})\Gamma(\mu+\alpha+1)}{\Gamma(\alpha)\Gamma(\mu+1)T^{\alpha+\mu} E_{\alpha,1+\frac{\beta}{\alpha},1+\frac{\beta+\mu}{\alpha}}(-\lambda_{k}T^{\alpha+\beta})}\right|^2
$$
$$
\leq C\sum\limits_{k=1}^n \frac{|\varphi_{k}|^2}{T^{\mu+1} |E_{\alpha, 1+\frac{\beta}{\alpha},1+\frac{\beta+\mu}{\alpha}}(-\lambda_k
T^{\alpha+\beta})|^2}.
$$
By virtue of (\ref{ML2}),
$$I_{2,n}\leq C\sum\limits_{k=1}^n \frac{\lambda_k^2 |\varphi_k|^2}{T^{2\mu+2}\left( \frac{1}{T^{\alpha+\beta}}\right)^2}
\leq CT^{2(\alpha+\beta+1-\mu)} \sum\limits_{k=1}^n \lambda_k^2  |\varphi_k|^2.
$$

Thus, since $\varphi, \Psi \in D(A)$, from estimates of $I_{1,
n},\, I_{2, n}$ and (\ref{Fj}) we finally have $f\in H$.

The fulfillment  of
Definition \ref{defIN} conditions   for function $u(t)$,
defined by the series (\ref{ux}) is proved in exactly the same
way as in Theorem \ref{th1}.

\textbf{Uniqueness.} It suffices to prove that the solution of the homogeneous inverse problem. By virtue of (2), (6):
\begin{equation}\label{eq100}
\partial_t^\alpha u(t) + t^{\beta}Au(t) = t^{\mu}f, \quad  t>0;
\end{equation}
\begin{equation}\label{in101}
 \lim\limits_{t\rightarrow +0}J_t^{\alpha-1} u(t) =0
\end{equation}
\begin{equation}\label{in111}
u(T) =0,
\end{equation}
therefore $u(t)\equiv 0$, and $f=0$.

Let  $u(t)$ be a solution to this problem and $u_k(t)=(u(t),
v_k)$. Then, by virtue of equation (\ref{eq100}) and the
selfadjointness of operator $A$,
$$
\partial_t^\alpha u_k(t)=(\partial_t^\alpha u(t), v_k)= -t^{\beta}(Au(t), v_k)
+(t^{\mu}f,v_k)=-t^{\beta}(u(t), Av_k) +(t^{\mu}f,v_k)
$$
\begin{equation}\label{eq102}
=-t^{\beta}(u(t), \lambda_k v_k) +t^{\mu}f_k=-\lambda_k t^{\beta}(u(t), v_k) +t^{\mu}f_k=-\lambda_k t^{\beta}u_k(t) +t^{\mu}f_k, t>0.
\end{equation}
Thus, taking into account (\ref{in101}),  we have the following
problem
$$
\partial_t^\alpha u_k(t) +\lambda_k t^{\beta}u_k(t)-t^{\mu}f_k=0,\quad t>0; \quad  J^{\alpha-1}_t u_k(+0) =0.
$$
A solution to this problem has the form (see,example, \cite{KST}, p. 246-247)
$$
u_k(t)=  \frac{\Gamma(\mu+1){f_{k}}}{\Gamma(\mu+\alpha+1)} t^{\alpha+\mu} E_{\alpha,1+\frac{\beta}{\alpha},1+\frac{\beta+\mu}{\alpha}}(-\lambda_{k}t^{\alpha+\beta}).
$$
Using (\ref{in111}),  we have
$$
u_k(T)=\frac{\Gamma(\mu+1){f_{k}}}{\Gamma(\mu+\alpha+1)} T^{\alpha+\mu} E_{\alpha,1+\frac{\beta}{\alpha},1+\frac{\beta+\mu}{\alpha}}(-\lambda_{k}T^{\alpha+\beta})=0.
$$
But, due to the properties of the Mittag-Leffler function,

$E_{\alpha,1+\frac{\beta}{\alpha},1+\frac{\beta+\mu}{\alpha}}(-\lambda_{k}T^{\alpha+\beta})\neq 0$. Therefore
$f_k=0$, for all $k\geq 1$.
In consequence, from the completeness
of the system of eigenfunctions $ \{v_k \} $, we finally obtain
$f=0$ and $u(t) \equiv 0$, as required.

\section{Conclusion}

Initial boundary value problems for subdiffusion equations have been studied by many mathematicians with different elliptic parts.  If the coefficients of the elliptic part depend on time, then either the variables do not separate, or the solution of the equation with respect to $t$ is not known.  In this regard, such equations are little studied, although various processes are modeled precisely by such equations.  The work of Yamamoto and Kubica \cite{Kubica} considers the initial-boundary value problem for the subdiffusion equation whose elliptic part depends on time.  The existence and uniqueness of a weak solution was proved by the method of a priori estimates.  In the present paper, for a certain class of equations whose coefficients depend on time, the existence and uniqueness of a generalized solution is proved by the Fourier method.
The paper also considers the inverse problem of determining the right side of the equation.

\par\bigskip\noindent
{\bf Acknowledgment.} The authors are grateful to  Alimov Sh.A.
for discussions of these results and acknowledge financial support from the Ministry of Innovative Development of the Republic of Uzbekistan, Grant No F-FA-2021-424.

\nocite{*}
\bibliography{aipsamp}

\end{document}